\documentclass[11pt,a4paper]{article}
\usepackage[T1]{fontenc}
\usepackage{lmodern}
\usepackage[margin=25mm]{geometry}
\usepackage{amsmath,amssymb,bm}
\usepackage{graphicx,booktabs,placeins}
\usepackage{microtype}
\usepackage[numbers,sort&compress]{natbib}
\usepackage[hidelinks]{hyperref}
\hypersetup{pdftitle={Phase Transitions and Modulation in Dense Oscillatory Associative Memories},pdfauthor={Shurui Liu and Zhuchun Li}}
\title{Phase Transitions and Modulation in Dense Oscillatory Associative Memories}
\author{Shurui Liu \quad Zhuchun Li\thanks{\href{mailto:lizhuchun@hit.edu.cn}{lizhuchun@hit.edu.cn}}\\[0.4em]
\normalsize Department of Mathematics, Harbin Institute of Technology,\\
\normalsize Harbin, 150001, China}
\date{}
\begin{document}
\maketitle
\begin{abstract}
We introduce high-order Winfree models for oscillatory associative memory. We characterize memory performance through the local-stability dichotomy
between stored and unstored patterns. Using this criterion, we show that crosstalk noise drives a cascade of ferromagnetic, mixed, and spin-glass transitions as the memory load increases. For both polynomial and exponential Winfree models, we derive the phase boundaries and develop auxiliary-field modulation schemes adapted to the different phases. We further present corresponding circuit designs and numerical simulations that support our theoretical predictions.
\end{abstract}

\section{Introduction}

Associative memory networks recover stored patterns from corrupted inputs. For a network of \(N\) units storing \(M\) binary patterns, the capacity scaling is conventionally characterized by the maximal growth of \(M\) with \(N\) for which corrupted inputs can still be reliably corrected to their corresponding stored patterns. The classical Hopfield model~\cite{doi:10.1073/pnas.79.8.2554} provides the foundational framework for attractor-based associative memory but has limited capacity scaling. To overcome this limitation, this framework has been extended to modern forms~\cite{krotov2016,demircigil2017,ramsauer2021hopfield}, broadly referred to as  Dense Associative Memories (DAMs)~\cite{krotov2016}. These models employ generalized energy functions with high-order interactions,
$
E=-\sum_{\mu=1}^{M}F\!\left(\sum_{k=1}^{N}\xi_k^\mu \sigma_k\right),
$
where $\xi_k^\mu=\pm1$ denotes the $k$th component of the $\mu$th stored binary pattern and $\sigma_k=\pm1 $ denote the $k$th neuron. Such generalized energies, together with the induced update rules, substantially enhance the corresponding capacity scaling. In particular, a polynomial choice $F(x)=x^n$ supports the scaling $M\sim N^{n-1}/\log N$~\cite{krotov2016}, while an exponential choice $F(x)=e^x$ yields the exponential scaling $M\sim c^N$~\cite{demircigil2017}.

In oscillatory neural networks, neurons are modeled as phase oscillators governed by coupled differential equations. Such continuous phase dynamics offer rich computational behavior~\cite{10.1063/1.5120412} and are naturally compatible with hardware implementations~\cite{TodriSanial2024}. To exploit these properties for associative memory, various frameworks have been proposed. Nishikawa \textit{et al.}\ applied the Kuramoto model~\cite{kuramoto1984}, showing that the introduction of a second-harmonic term stabilizes otherwise fragile memories and supports the capacity scaling $M(N)\sim N/\log N$~\cite{PhysRevLett.92.108101}. More recently, N\"{a}gerl and Berloff introduced a generalized Kuramoto model incorporating $2l$-body interactions~\cite{NagerlBerloff2025}. By analyzing the specific case of simultaneous pairwise and quartic couplings, they identified a tricritical point in the phase diagram and showed that such high-order terms can produce superlinear capacity scaling. In parallel, Li and Zhao established an associative-memory framework based on Winfree-type oscillators in the pairwise interaction regime~\cite{10.1098/rspa.2025.0154}.

In this work, we develop a Winfree-based framework for dense oscillatory associative memories with high-order interactions.
We derive the dense Winfree models from the gradient flow of a real-valued energy $E = -\sum_{\mu=1}^{M} F(\sum_{k=1}^{N}\xi_{k}^{\mu}\cos\theta_k)$, which energy parallels the DAMs form . By specifying $F(x)$, we obtain two model classes: the \textit{Polynomial Winfree model} with $F(x) = \frac{x^m}{mN^{m-1}}$ (classified into odd and even cases), and the \textit{Exponential Winfree model} with $F(x) = e^{x-N}$.

The energy landscape admits every stored or unstored binary pattern as a
stationary state, with phases \(0\) and \(\pi\) encoding \(+1\) and \(-1\),
respectively. These states are generally not distinguished well by their local stability. We therefore introduce a conservative baseline auxiliary field $S_p(\epsilon)=-\epsilon(N)\sin\theta_p\cos\theta_p$ to render stored patterns stable and unstored patterns unstable, while preserving the gradient-flow structure:
$
\dot{\theta}_p=-\frac{\partial E}{\partial \theta_p}+S_p(\epsilon).
$
The key quantity is therefore the required auxiliary-field parameter $\epsilon(N)$, whose scaling reveals a cascade of phase transitions in the oscillator network as the memory load $M$ increases: a ferromagnetic phase (FM), in which the required parameter remains infinitesimal; a mixed phase (MP), in which the required parameter diverges with $N$; and a spin-glass phase (SG), in which no admissible choice of $\epsilon(N)$ exists.

We derive the phase boundaries of these transitions, summarized in Table~\ref{tab:phase_boundaries}.
These phase boundaries, set by the local-stability dichotomy between stored and unstored patterns, mark qualitative changes in the memory performance of the oscillator network. For practical reason, we replace the baseline auxiliary field \(S_p(\epsilon)\) with tailored fields modulates the energy landscape and enlarges the attraction basins of stored patterns. We further provide circuit implementations of the Winfree network and its auxiliary fields, while simulations highlight the practical potential of memory operation in both the FM and MP.

\begin{table}[tbp]\centering
\caption{Asymptotic phase boundaries \(M(N)\)}
\label{tab:phase_boundaries}
\renewcommand{\arraystretch}{1.6}
\begin{tabular}{lcc}\toprule
Model & FM$\to$MP & MP$\to$SG \\\midrule
Even Winfree &
$\frac{N^{m/2}}{(m-1)!!}$ &
$\frac{12N^{m-1}\ln N}{(2m-3)!!\pi^2}$ \\
Odd Winfree &
$\frac{N^{m-1}}{2(2m-3)!!\ln N}$ &
$\frac{12N^{m-1}\ln N}{(2m-3)!!\pi^2}$ \\
Exponential Winfree &
$\sim 1.762^N$ &
$\sim 1.964^N$ \\
\bottomrule\end{tabular}

\end{table}

\section{Phase-transition boundaries}
We first clarify the physical mechanism underlying the phase transitions in dense
oscillatory associative-memory networks. Let
\(\boldsymbol{\xi}_N=[\xi^1,\ldots,\xi^M]\in
\{\pm1\}^{N\times M}\) denote the random memory task, whose entries
\(\xi_k^\mu\) are independent Rademacher variables. A load scaling
\(M=M(N)\) is said to be memory-admissible if the network exhibits
\textit{perfect memory}. For oscillatory associative-memory networks,
perfect memory means that, as \(N\to\infty\), the proportion of locally
stable patterns among the \(M\) stored patterns and the proportion of locally
unstable patterns among the unstored patterns both converge to
one in probability over the memory task. We characterize memory performance
in terms of this local-stability dichotomy, which is encoded by the Jacobian spectra
at the corresponding equilibria.

Consider a randomly selected stored pattern $\xi^\upsilon$ and an unstored pattern $\eta$. Let $\lambda_{\max}(\mathbf J(\boldsymbol{\xi}_N)|_{\xi^\upsilon})$ and $\lambda_{\max}(\mathbf J(\boldsymbol{\xi}_N)|_{\eta})$ denote the largest eigenvalues of the Jacobian of the gradient flow, $\dot{\theta}_p=-\frac{\partial E}{\partial \theta_p},$ evaluated at the corresponding equilibria. In the framework of dense Winfree network, the Jacobian is diagonal at any binary equilibrium. The baseline auxiliary field $S_p(\epsilon)$ shifts the entire Jacobian spectrum at any binary equilibrium by $\epsilon(N)$ toward the negative real axis. A stored pattern $\xi^{\upsilon}$ is therefore stabilized with probability $\mathbb P\!\left(\lambda_{\max}(\mathbf J(\boldsymbol{\xi}_N)|_{\xi^\upsilon})<\epsilon(N)\right),$ whereas an unstored pattern $\eta$ is destabilized with probability $\mathbb P\!\left(\lambda_{\max}(\mathbf J(\boldsymbol{\xi}_N)|_{\eta})>\epsilon(N)\right).$ We define the overall memory discriminability, $\mathbb D$, as the product of these two probabilities. Perfect memory requires \(\mathbb D\to1\) as \(N\to\infty\).

For a stored pattern \(\xi^\upsilon\), the contribution of the memory
task to the local stability separates into the deterministic signal
generated by \(\xi^\upsilon\) itself and the crosstalk noise generated
by the remaining \(M-1\) patterns
\(\{\xi^\mu\}_{\mu\ne\upsilon}\). At an unstored pattern \(\eta\),
there is no matching memory to provide a signal, and all \(M\) stored patterns instead
generate the background noise. Under the corresponding binary gauge transformations, stored and
unstored binary patterns can each be mapped to reference state
\(\mathbf 1_N\) without changing their Jacobian spectra, while the
noise-generating patterns in the random memory task
\(\boldsymbol{\xi}_N\) remain i.i.d.\ Rademacher. Consequently, the contributors of the noise at stored and unstored patterns
differ only by a single pattern, which is negligible as
\(N\to\infty\); the background noise and the crosstalk noise are therefore regarded as asymptotically equivalent. We define
\(\mathbf J:=\left.\mathbf J(\boldsymbol{\xi}_N)\right|_{\mathbf 1_N}\)
as the reference random Jacobian and refer to its largest eigenvalue
\(\lambda_{\max}(\mathbf J)\) as the crosstalk-noise statistic. We
denote the mean and variance of this statistic by \(\mu_\lambda\) and
\(\sigma_\lambda^2\), respectively. For the polynomial and exponential Winfree models studied here, the
signal term shifts the Jacobian spectrum by \(-1\) relative to the reference Jacobian spectrum, so the stability
probability of a stored pattern is asymptotically
\(\mathbb P\!\left(\lambda_{\max}(\mathbf J)-1<\epsilon(N)\right)\).

This spectral characterization reveals the physical mechanism of the
two phase transitions. In the thermodynamic limit (\(N\to\infty\)),
both the mean and variance of \(\lambda_{\max}(\mathbf J)\) remain negligible in
the FM, so an infinitesimal \(\epsilon(N)\) suffices to preserve
perfect memory.  As the memory load \(M\) increases, the
FM--MP transition is driven by the mean \(\mu_\lambda\), whose
magnitude becomes macroscopic. We term this the mean-driven transition
and locate the FM--MP boundary by the condition
\(|\mu_\lambda|=1\). Throughout the
MP, \(|\mu_\lambda|\) diverges while \(\sigma_\lambda^2\) remains
negligible, leaving a spectral gap of width \(1\) between the largest
eigenvalues associated with unstored and stored patterns. A
diverging auxiliary-field parameter
\(\epsilon(N)=\mu_\lambda-\tfrac12\), centered within this gap,
therefore preserves perfect memory.

The subsequent MP--SG transition is driven by the variance
\(\sigma_\lambda^2\), which grows from negligible to macroscopic scale. We term
this the variance-driven transition and locate the MP--SG boundary by
the condition \(\sigma_\lambda^2=1\). In the SG phase, the spectral
distributions become sufficiently broad so that no choice of the
auxiliary-field parameter \(\epsilon(N)\) can reliably distinguish
the largest eigenvalues associated with stored and unstored patterns,
making perfect memory impossible.

We note that, since \(M\ll 2^N\), the proportion of locally unstable
patterns among the unstored patterns may be replaced by that among all
binary patterns, with an error of at most \(M/2^N=o(1)\), thereby
simplifying the analysis.

\begin{figure}[tbp]
\centering
\includegraphics[width=0.65\linewidth]{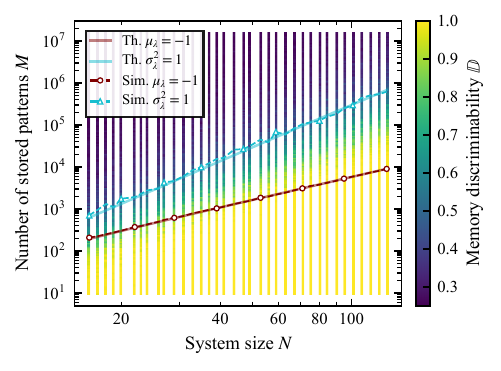}
\caption{
Numerical phase diagram for the even polynomial Winfree model
(\(m=4\)). The background shows the experimental memory
discriminability \(\mathbb D\). Solid curves are the theoretical
boundaries \(\mu=-1\) and \(\sigma^2=1\); dashed curves with markers are
the corresponding experimental boundaries.
}
\label{fig:phase_diagram_winfree_even}
\end{figure}

Hence, the two phase transitions are dictated, respectively, by the mean and
variance of the crosstalk-noise statistic
\(\lambda_{\max}(\mathbf J)\). We now derive the phase boundaries for each model.
We begin with the polynomial Winfree model,
\begin{equation}
\dot{\theta}_p
=
-\frac{\sin\theta_p}{N^{m-1}}
\sum_{j_1,\ldots,j_{m-1}=1}^{N}
C_{j_1\cdots j_{m-1}p}
\cos\theta_{j_1}\cdots\cos\theta_{j_{m-1}},
\label{eq:poly_winfree_main}
\end{equation}
where
\(C_{j_1\cdots j_{m-1}p}
=
\sum_{\mu=1}^{M}
\xi_{j_1}^{\mu}\cdots
\xi_{j_{m-1}}^{\mu}\xi_p^\mu\)
is defined by the high-order Hebbian rule. At the equilibrium associated with
\(\eta\in\{\pm1\}^N\), the Jacobian is diagonal and satisfies
\(J_{pq}|_{\eta}
=
-\frac{\delta_{pq}}{N^{m-1}}
\sum_{\mu=1}^{M}
\langle\xi^\mu,\eta\rangle^{m-1}\xi_p^\mu\eta_p\). Under the binary gauge transformation, the reference random Jacobian becomes $J_{pq} = -\frac{\delta_{pq}}{N^{m-1}}\sum_{\mu=1}^{M}
\langle \xi^\mu,\mathbf 1_N\rangle^{m-1}\xi_p^\mu.$ Its diagonal entries are
\(\Gamma_j=\sum_{\mu=1}^{M}G_{\mu j}\), with
\(G_{\mu j}
=
-\frac{1}{N^{m-1}}
(\xi_j^\mu+\sum_{k\ne j}\xi_k^\mu)^{m-1}\xi_j^\mu\).
Hence the crosstalk-noise statistic is
\(\lambda_{\max}(\mathbf J)=\max_{1\le j\le N}\Gamma_j\).

For even \(m\ge4\), a binomial expansion gives the mean and variance of
\(G_{\mu j}\), which do not depend on \(\mu\) or \(j\):
\(\mu_G\sim\frac{-(m-1)!!}{N^{m/2}}\) and
\(\sigma_G^2\sim\frac{(2m-3)!!}{N^{m-1}}\), respectively.
Since $\Gamma_j$ is a sum of $M$ independent terms, the central limit theorem (CLT) implies that each $\Gamma_j$ is asymptotically Gaussian, with mean $\mu_\Gamma=M\mu_G
$ and variance $\sigma_\Gamma^2=M\sigma_G^2
$. Moreover, for all $j\neq l$, the correlation coefficient between $\Gamma_j$ and $\Gamma_l$ is $\rho_\Gamma \sim \frac{1}{N}\!\left[2(m-1)-\frac{[(m-1)!!]^2}{(2m-3)!!}\right],$ so that $\rho_\Gamma$ vanishes in the thermodynamic limit.

Therefore, we replace all the \(\Gamma_j\) by \(N\) independent
Gaussian variables with mean \(\mu_\Gamma\) and variance
\(\sigma_\Gamma^2\). Extreme-value theory then gives an asymptotic Gumbel distribution
for \(\lambda_{\max}(\mathbf J)\), with location parameter
\(b=\mu_\Gamma+\sigma_\Gamma b_N\) and scale parameter
\(a=\sigma_\Gamma/\sqrt{2\ln N}\), where
\(b_N\sim\sqrt{2\ln N}\). Moreover, the mean and variance of the crosstalk-noise statistic are
\(\mu_\lambda=b+a\gamma\) and
\(\sigma_\lambda^2=\pi^2a^2/6\), respectively, where
\(\gamma\) is the Euler--Mascheroni constant. The mean-driven FM--MP boundary, determined by
\(|\mu_\lambda|=1\), therefore is at
\(M\sim N^{m/2}/(m-1)!!\). The variance-driven MP--SG boundary,
determined by \(\sigma_\lambda^2=1\), is likewise at
\(M\sim 12N^{m-1}\ln N/[(2m-3)!!\pi^2]\). Figure~\ref{fig:phase_diagram_winfree_even} shows the numerical phase
diagram for the even polynomial Winfree model.

For odd \(m\ge3\), analogous calculations give
\(\mu_G=0\) and
\(\sigma_G^2\sim(2m-3)!!/N^{m-1}\), hence
\(\mu_\Gamma=0\), \(\sigma_\Gamma^2=M\sigma_G^2\), and each
\(\Gamma_j\) is asymptotically Gaussian by the CLT. Since
\(\rho_\Gamma\sim2(m-1)/N\to0\) for \(j\ne l\), the same Gaussian
replacement and extreme-value argument yields an asymptotic Gumbel
distribution for \(\lambda_{\max}(\mathbf J)\). The mean-driven FM--MP
and variance-driven MP--SG boundaries therefore occur at
\(M\sim N^{m-1}/[2(2m-3)!!\ln N]\) and
\(M\sim12N^{m-1}\ln N/[(2m-3)!!\pi^2]\), respectively.

We next consider the exponential Winfree model,
\begin{equation}
\dot{\theta}_p
=
-e^{-N}\sin\theta_p
\sum_{\mu=1}^{M}
\xi_p^\mu
\exp\!\left(
\sum_{k=1}^{N}\xi_k^\mu\cos\theta_k
\right).
\label{eq:exp_winfree_main}
\end{equation}
At a binary equilibrium \(\eta\in\{\pm1\}^N\), the Jacobian is diagonal,
with

\[
J_{pq}|_{\eta}
=
-\delta_{pq}
\sum_{\mu=1}^{M}
\exp\left( -N+\sum_k\xi_k^\mu\eta_k\right) \xi_p^\mu\eta_p.
\]

After the binary gauge transformation, the reference Jacobian becomes

\[
J_{pq}
=
-\delta_{pq}
\sum_{\mu=1}^{M}
\exp\left( -N+\sum_k\xi_k^\mu\right) \xi_p^\mu.
\]

For simplicity, we use the same notation for the polynomial and exponential models. The crosstalk-noise statistic is \(\lambda_{\max}(\mathbf J)=\max_j\Gamma_j\), where \(\Gamma_j=\sum_{\mu}G_{\mu j}\) and
$
G_{\mu j} = -\exp\left( -N+\sum_k\xi_k^\mu\right) \xi_j^\mu.
$

Direct averaging over the independent Rademacher variables gives
\[
\mu_G=-e^{-N}(\cosh 1)^{N-1}\sinh 1,
\qquad \sigma_G^2\sim e^{-2N}(\cosh 2)^N.
\] Independence across the memory
index \(\mu\) then gives
\(\mu_\Gamma=M\mu_G\) and
\(\sigma_\Gamma^2=M\sigma_G^2\).
Moreover, for \(j\ne l\), the correlation coefficient between
\(\Gamma_j\) and \(\Gamma_l\) satisfies
\(\rho_\Gamma\sim\tanh^2 2\). The Berry--Esseen theorem implies that
\(\Gamma_j\) is asymptotically Gaussian for
\(M=c^{N+o(N)}\) with \(c>c_{\mathrm{BE}}\approx1.903\).

We now determine the base \(c_{\mathrm{FM\text{-}MP}}\) of the exponential load \(M=c^{N+o(N)}\) at the mean-driven boundary \(|\mu_\lambda|=1\). At exponential order, setting
\(|\mu_\Gamma|=1\) yields a candidate base
\(c=e/\cosh 1\). At this base,
$$0\le\mu_\lambda-\mu_\Gamma
=
\mathbb E[\max_{1\le j\le N}(\Gamma_j-\mu_\Gamma)]
\le
\mathbb E[(\sum_{j=1}^{N}(\Gamma_j-\mu_\Gamma)^2)^{1/2}]
\le
\sqrt{N(\sigma_\Gamma^2 +\mu_\Gamma^2)}
=
o(1).$$
The first upper bound follows from
\(\max_j x_j\le(\sum_jx_j^2)^{1/2}\), while the second follows from
the Cauchy--Schwarz inequality
\(\mathbb E[X]\le\sqrt{\mathbb E[X^2]}\) for any nonnegative
square-integrable random variable \(X\), applied to
\(X=(\sum_{j=1}^{N}(\Gamma_j-\mu_\Gamma)^2)^{1/2}\).
The final \(o(1)\) follows from \(\cosh 2<e\cosh 1\).
Thus, \(\mu_\lambda=\mu_\Gamma+o(1)\), so the condition
\(|\mu_\lambda|=1\) confirms this candidate, yielding
\(c_{\mathrm{FM\text{-}MP}}\approx1.762\).

Since correlation coefficient \(\rho_\Gamma\) remains macroscopic in the thermodynamic limit, we use
the equicorrelated Gaussian representation
\(\Gamma_j=\mu_\Gamma+\sigma_{\mathrm{com}}Z_0+\sigma_{\mathrm{res}}Z_j\), where
\(Z_0,Z_1,\ldots,Z_N\) are independent standard normal variables. Here \(Z_0\) represents the common fluctuation shared by all diagonal entries, while the \(Z_j\)'s represent the residual independent fluctuations, with \(\sigma_{\mathrm{com}}^2=\rho_\Gamma\,\sigma_\Gamma^2\) and \(\sigma_{\mathrm{res}}^2=(1-\rho_\Gamma)\sigma_\Gamma^2\). Hence \(\lambda_{\max}(\mathbf J)= \mu_\Gamma+\sigma_{\mathrm{com}}Z_0 + \sigma_{\mathrm{res}}\max_j Z_j\). The residual extreme-value contribution has variance \(\sigma_{\mathrm{res}}^2\pi^2/(12\ln N)\), which is smaller than the common contribution \(\sigma_{\mathrm{com}}^2=\rho_\Gamma\,\sigma_\Gamma^2\) by a factor of order \(1/\ln N\). Therefore, $\sigma_\lambda^2$ becomes \(O(1)\) when \(\sigma_\Gamma^2=O(1)\), giving \(c_{\mathrm{MP\text{-}SG}}=e^2/\cosh 2\approx1.964\). Since \(c_{\mathrm{MP\text{-}SG}}>c_{\mathrm{BE}} \approx 1.903\), the Gaussian
approximation is self-consistent at the variance-driven boundary.

\section{Modulation}
The phase boundaries characterize the local-stability contrast between stored
and unstored binary patterns. Although the baseline auxiliary field
can exploit this contrast to stabilize the stored patterns while destabilizing
the unstored patterns, local stability alone does not guarantee a practically wide
basin of attraction. We call the corrupted input the query and the stored
pattern to be recovered its target.
According to the different energy landscapes of polynomial and exponential
models in the FM and MP, We tailor auxiliary
fields to obtain large restoring regions of target and denote the corresponding auxiliary energy as $E^{\rm aux}$. We next explain the notion of restoring regions.

We gauge a prescribed stored target pattern to \(\mathbf1_N\) and
consider a prescribed query \(\mathbf y^0\), with \(y_p^0=\cos\theta_p^0\).
For each query coordinate
\(y_p^0<1\), let \(y_p\) increase from \(y_p^0\) to the target value \(1\),
keeping all other query coordinates unchanged.
Write the total energy along this interval as \(E(y_p)\).
The \(p\)th coordinate is counted as failed if \(dE(y_p)/dy_p\ge0\)
at some point in \([y_p^0,1]\). We say a prescribed query is in the restoring region of its target if
the expected number of failed coordinates and the probability that
at least one failed coordinate exists tend to zero as \(N\to\infty\).

At a prescribed unstored binary pattern \(\eta\), we again vary the
\(p\)th coordinate and fix all others. The coordinate is counted as
escaping if the energy strictly decreases as \(y_p\) moves away from
\(\eta_p\) along an interval of nonzero length independent of \(N\).
For each fixed memory task, we sample an unstored binary pattern
uniformly. We say that the network satisfies the criterion for
suppressing spurious binary memories if the expected number of
escaping coordinates diverges and the probability that at least one
escaping coordinate exists tends to one as \(N\to\infty\), both in
probability over memory tasks.

We first consider the polynomial Winfree model with \(m>2\) in the
FM. For any fixed \(\chi\in\left( 0,\,1/2\right) \), choose the auxiliary energy as
\(E^{\rm quad}
=\frac{(1-2\chi)^{m-1}}{4}\sum_{j=1}^{N}\cos^2\theta_j\).
For loads asymptotically below the FM--MP scale, the restoring region
around the gauged target \(\mathbf1_N\) contains
\(\mathcal{S}^0 = \{\mathbf y^0\in[-1,1]^N:h(\mathbf y^0)\ge h_0:=1-2\chi>0\}\)
, where
\(h(\mathbf y^0)=N^{-1}\sum_j y_j^0\) is the overlap with the target.
Binary queries with at most \(\chi N\) phase flips relative to the
target are in $\mathcal{S}^0$.

For the odd polynomial Winfree model in the MP, we use the tailored
auxiliary energy \(E^{\rm aux}=E^{\rm abs}\), with absolute-value energy
$
E^{\rm abs}
=-\left(g_N-\frac12\right)\sum_{j=1}^{N}|\cos\theta_j|.
$
Here \(v_N=[M(2m-3)!!/N^{m-1}]^{1/2}\) and \(g_N=b_Nv_N\),
with \(b_N\sim\sqrt{2\ln N}\) as above. The restoring region contains every prescribed query satisfying
\(0<y_p^0\le1\) for all \(p\). Thus, the field supports recovery
of analog inputs whose perturbations preserve the target signs.
In contrast, the target is generally unstable under binary
coordinate updates in the corresponding polynomial DAM, even when
the binarized query coincides with the target.

Our polynomial energy allows repeated oscillator indices, as in the
higher-order Kuramoto model~\cite{NagerlBerloff2025}. This differs
from the interactions among distinct neurons studied in
Refs.~\cite{BaldiVenkatesh1987,Gardner1987,AbbottArian1987}. For odd $m$, the terms with repeated indices naturally cancel out. For even \(m\), the expansion of the single-pattern energy
\(-(\sum_k\xi_k^\mu y_k)^m/(mN^{m-1})\) contains terms in which
every index occurs an even number of times.
At \(m=4\), for example, the indices \((i,i,j,j)\), with \(i\ne j\),
give a term proportional to \(-y_i^2y_j^2\).
Summing all such terms over the \(M\) memories gives a deterministic
self-energy minimized at every binary pattern. We cancel this
contribution by adding
$
E^{\rm self}(\mathbf y)=\frac{M}{mN^{m-1}}\Psi_m(\mathbf y),
$
where \(\Psi_m(\mathbf y)=\mathbb E_{\xi}[(\sum_k\xi_k y_k)^m]\)
for one independent Rademacher pattern \(\xi\).
For $M\ll N^{m-1}/\ln N $, the compensated
landscape is FM-like. With \(E^{\rm quad}\) as in the FM, choosing
\(E^{\rm aux}=E^{\rm self}+E^{\rm quad}\) supports recovery from
the same positive-overlap region, including binary queries with a
fixed fraction of phase flips.
For higher loads
\(N^{m-1}/\ln N\ll M\ll N^{m-1}\ln N\),
\(E^{\rm aux}=E^{\rm self}+E^{\rm abs}\), with \(E^{\rm abs}\)
as in the odd polynomial Winfree model, instead supports recovery
of prescribed queries satisfying \(0<y_p^0\le1\) for all \(p\).

For the exponential model, choose
\(E^{\rm aux}=E^{\rm self}+E^{\rm quad}\), where
\(E^{\rm self}(\mathbf y)=Me^{-N}\prod_k\cosh y_k\) cancels
the self-energy contribution and
\(E^{\rm quad}=\frac14e^{-2\chi N}\sum_j\cos^2\theta_j\).
For \(M=c^{N+o(N)}\), with fixed \(1<c<e^2/\cosh2\) and
\(0\le\chi<\frac14\ln[e^2/(c\cosh2)]\), the restoring region
contains every prescribed query satisfying \(h(\mathbf y^0)\ge1-2\chi\).
The network therefore supports coordinate recovery toward the target
for prescribed binary queries with at most \(\chi N\) phase flips,
as well as continuous queries satisfying the same overlap condition.
Every fixed exponential base \(c\) inside either the FM or MP thus
allows a strictly positive phase-flip fraction \(\chi\).
For the choice \(\chi=\chi_N\to0\), the quadratic field alone
suffices in the FM,
whereas the MP construction requires additional self-energy
compensation for the stated recovery guarantee.

\section{Circuit implementation}
Dense Winfree models with a general differentiable interaction
function \(F\) share a common circuit architecture.
In this circuit section, \(x_p=\cos\theta_p\) and
\(y_p=\sin\theta_p\) are normalized state voltages. With
\(S_\mu=\sum_j\xi_j^\mu x_j\), the energy
\(E=-\sum_\mu F(S_\mu)+E^{\rm aux}(\mathbf x)\) generates
\begin{equation}
\dot\theta_p=-\sin\theta_p A_p
-\frac{\partial E^{\rm aux}}{\partial\theta_p},
\label{eq:general_winfree_modulation}
\end{equation}
where \(A_p=\sum_\mu\xi_p^\mu F'(S_\mu)\).

\begin{figure}[tbp]
\centering
\includegraphics[width=0.60\linewidth]{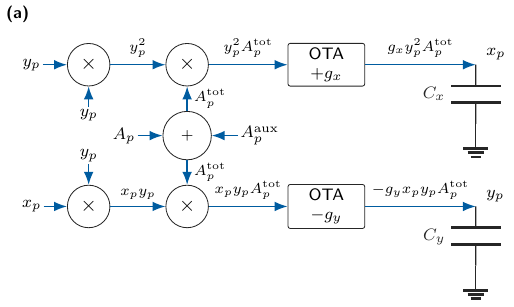}
\vspace{1mm}
\includegraphics[width=0.60\linewidth]{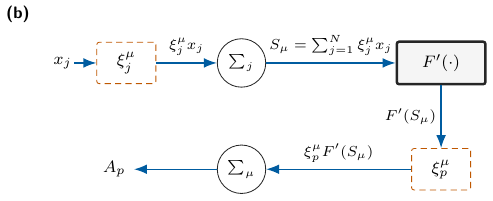}
\caption{Analog architecture for dense Winfree dynamics.
(a) Two-capacitor state cell driven by
\(A_p^{\rm tot}=A_p+A_p^{\rm aux}\).
OTA denotes an operational transconductance amplifier.
(b) Signed summation and the nonlinear response \(F'\) generate
the interaction field from shared memory overlaps.}
\label{fig:winfree_circuit}
\end{figure}

For a fixed memory task, the circuit in Fig.~\ref{fig:winfree_circuit}(b)
implements different dense Winfree models by replacing the \(F'\) blocks.
This simple construction motivates retaining the full \(F(S_\mu)\),
including any self-energy contribution, in the base model.

The two-capacitor representation avoids explicit sine or cosine
evaluation. Setting \(A_p^{\rm aux}=-\partial E^{\rm aux}/\partial x_p\)
and
\(A_p^{\rm tot}=A_p+A_p^{\rm aux}\), the dynamics
becomes \(\dot x_p=y_p^2A_p^{\rm tot}\) and
\(\dot y_p=-x_py_pA_p^{\rm tot}\), preserving
\(x_p^2+y_p^2=1\) for initial states on the unit circle.
Multipliers and two capacitor integrators with matched time scales
implement these equations directly [Fig.~\ref{fig:winfree_circuit}(a)].
The chosen auxiliary-field module supplies
\(A_p^{\rm aux}\) through the same field summing node.

\begin{figure}[tbp]
\centering
\includegraphics[width=0.6\linewidth]{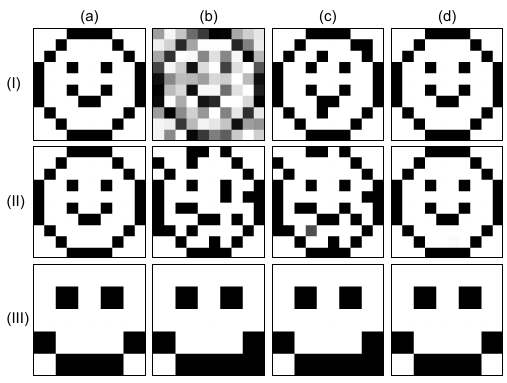}
\caption{Rows: (I) \(m=3\), \(N=100\), \(M=10^4\);
(II) \(m=4\), \(N=100\), \(M=7168\);
(III) exponential, \(N=25\), \(M=3\times10^5\).
Columns: (a) target, (b) query, (c) output with the baseline auxiliary field,
and (d) output with the tailored auxiliary field.}
\label{fig:three_model_circuit_retrieval}
\end{figure}

Figure~\ref{fig:three_model_circuit_retrieval} illustrates retrieval
of a smiley face target in ideal circuit simulations. The other
stored patterns are randomly generated.
For each model, the quadratic baseline and tailored auxiliary field
use the same memory library and initial query.
Grayscale denotes \(x_p=\cos\theta_p\).
Binary queries in rows (II)--(III) are shown before the small
initialization phase offsets.
The absolute-value field retrieves the target from a continuously
perturbed, sign-preserving query in the odd model [row (I)].
Self-energy compensation combined with quadratic modulation corrects
phase-flipped queries in the even and exponential models [rows (II)--(III)].
In all three examples, the tailored fields recover every target
sign, whereas the quadratic baselines leave residual errors.

\FloatBarrier

\begingroup\small
\begin{thebibliography}{13}
\providecommand{\natexlab}[1]{#1}
\providecommand{\url}[1]{\texttt{#1}}
\expandafter\ifx\csname urlstyle\endcsname\relax
  \providecommand{\doi}[1]{doi: #1}\else
  \providecommand{\doi}{doi: \begingroup \urlstyle{rm}\Url}\fi

\bibitem[Hopfield(1982)]{doi:10.1073/pnas.79.8.2554}
J.~J. Hopfield.
\newblock Neural networks and physical systems with emergent collective
  computational abilities.
\newblock \emph{Proceedings of the National Academy of Sciences}, 79\penalty0
  (8):\penalty0 2554--2558, 1982.
\newblock \doi{10.1073/pnas.79.8.2554}.
\newblock URL \url{https://www.pnas.org/doi/abs/10.1073/pnas.79.8.2554}.

\bibitem[Krotov and Hopfield(2016)]{krotov2016}
D.~Krotov and J.~J. Hopfield.
\newblock Dense associative memory for pattern recognition.
\newblock In D.~Lee, M.~Sugiyama, U.~Luxburg, I.~Guyon, and R.~Garnett,
  editors, \emph{Advances in Neural Information Processing Systems}, volume~29,
  pages 1172--1180. Curran Associates, Inc., 2016.
\newblock URL
  \url{https://proceedings.neurips.cc/paper_files/paper/2016/file/eaae339c4d89fc102edd9dbdb6a28915-Paper.pdf}.

\bibitem[Demircigil et~al.(2017)Demircigil, Heusel, L{\"o}we, Upgang, and
  Vermet]{demircigil2017}
M.~Demircigil, J.~Heusel, M.~L{\"o}we, S.~Upgang, and F.~Vermet.
\newblock On a model of associative memory with huge storage capacity.
\newblock \emph{Journal of Statistical Physics}, 168\penalty0 (2):\penalty0
  288--299, 2017.
\newblock \doi{10.1007/s10955-017-1806-y}.
\newblock URL \url{https://doi.org/10.1007/s10955-017-1806-y}.

\bibitem[Ramsauer et~al.(2021)Ramsauer, Sch{\"a}fl, Lehner, Seidl, Widrich,
  Adler, Gruber, Holzleitner, Kreil, Kopp, Klambauer, Brandstetter, and
  Hochreiter]{ramsauer2021hopfield}
H.~Ramsauer, B.~Sch{\"a}fl, J.~Lehner, P.~Seidl, M.~Widrich, T.~Adler,
  L.~Gruber, M.~Holzleitner, D.~Kreil, M.~Kopp, G.~Klambauer, J.~Brandstetter,
  and S.~Hochreiter.
\newblock Hopfield networks is all you need.
\newblock In \emph{International Conference on Learning Representations
  (ICLR)}. OpenReview.net, 2021.
\newblock URL \url{https://openreview.net/forum?id=tL89RnzIiCd}.

\bibitem[Csaba and Porod(2020)]{10.1063/1.5120412}
G.~Csaba and W.~Porod.
\newblock Coupled oscillators for computing: A review and perspective.
\newblock \emph{Applied Physics Reviews}, 7\penalty0 (1):\penalty0 011302, 01
  2020.
\newblock ISSN 1931-9401.
\newblock \doi{10.1063/1.5120412}.
\newblock URL \url{https://doi.org/10.1063/1.5120412}.

\bibitem[Todri-Sanial et~al.(2024)Todri-Sanial, Delacour, Abernot, and
  Sabo]{TodriSanial2024}
A.~Todri-Sanial, C.~Delacour, M.~Abernot, and F.~Sabo.
\newblock Computing with oscillators from theoretical underpinnings to
  applications and demonstrators.
\newblock \emph{npj Unconventional Computing}, 1:\penalty0 14, 2024.
\newblock \doi{10.1038/s44335-024-00015-z}.

\bibitem[Kuramoto(1984)]{kuramoto1984}
Y.~Kuramoto.
\newblock \emph{Chemical Oscillations, Waves, and Turbulence}.
\newblock Springer-Verlag, Berlin, 1984.

\bibitem[Nishikawa et~al.(2004)Nishikawa, Lai, and
  Hoppensteadt]{PhysRevLett.92.108101}
T.~Nishikawa, Y.-C. Lai, and F.~C. Hoppensteadt.
\newblock Capacity of oscillatory associative-memory networks with error-free
  retrieval.
\newblock \emph{Phys. Rev. Lett.}, 92\penalty0 (10):\penalty0 108101, Mar 2004.
\newblock \doi{10.1103/PhysRevLett.92.108101}.
\newblock URL \url{https://link.aps.org/doi/10.1103/PhysRevLett.92.108101}.

\bibitem[Nagerl and Berloff(2025)]{NagerlBerloff2025}
J.~Nagerl and N.~G. Berloff.
\newblock Higher-order {Kuramoto} oscillator network for dense associative
  memory, 2025.
\newblock URL \url{https://arxiv.org/abs/2507.21984}.

\bibitem[Li and Zhao(2025)]{10.1098/rspa.2025.0154}
Z.~Li and X.~Zhao.
\newblock Associative-memory network of {Winfree}-type oscillators and binary
  pattern retrieval: theory and algorithm.
\newblock \emph{Proceedings of the Royal Society A: Mathematical, Physical and
  Engineering Sciences}, 481\penalty0 (2320):\penalty0 20250154, 08 2025.
\newblock ISSN 1364-5021.
\newblock \doi{10.1098/rspa.2025.0154}.
\newblock URL \url{https://doi.org/10.1098/rspa.2025.0154}.

\bibitem[Baldi and Venkatesh(1987)]{BaldiVenkatesh1987}
P.~Baldi and S.~S. Venkatesh.
\newblock Number of stable points for spin-glasses and neural networks of
  higher orders.
\newblock \emph{Phys. Rev. Lett.}, 58\penalty0 (9):\penalty0 913--916, 1987.
\newblock \doi{10.1103/PhysRevLett.58.913}.

\bibitem[Gardner(1987)]{Gardner1987}
E.~Gardner.
\newblock Multiconnected neural network models.
\newblock \emph{J. Phys. A: Math. Gen.}, 20\penalty0 (11):\penalty0 3453--3464,
  1987.
\newblock \doi{10.1088/0305-4470/20/11/046}.

\bibitem[Abbott and Arian(1987)]{AbbottArian1987}
L.~F. Abbott and Y.~Arian.
\newblock Storage capacity of generalized networks.
\newblock \emph{Phys. Rev. A}, 36\penalty0 (10):\penalty0 5091--5094, 1987.
\newblock \doi{10.1103/PhysRevA.36.5091}.

\end{thebibliography}

\endgroup
\end{document}